\documentclass[runningheads]{llncs}
\usepackage{amsmath} 
\usepackage{amssymb} 
\usepackage{graphicx}
\usepackage[ruled]{algorithm2e}
\usepackage{hyperref}
\usepackage{xcolor}
\newcommand{\rstar}{[0,\infty]}
\newcommand{\order}{\mathcal{O}}
\DeclareMathOperator*{\argmin}{argmin}
\DeclareMathOperator*{\argmax}{argmax}

\begin{document}
\title{Sparse Nerves in Practice}
%
%
\author{Nello Blaser \inst{1} 
\and
Morten Brun\inst{1}}
\authorrunning{N. Blaser and M. Brun}
%
\institute{Department of Mathematics, University of Bergen,
  All\'egaten 41, Bergen, Norway}
\maketitle              
\begin{abstract}

  Topological data analysis combines machine learning with methods
  from algebraic topology. Persistent homology, a method to
  characterize topological features occurring in data at multiple
  scales is of particular interest. A major obstacle to the
  wide-spread use of persistent homology is its computational
  complexity. In order to be able to calculate persistent homology of
  large datasets, a number of approximations can be applied in order
  to reduce its complexity.  We propose algorithms for calculation of
  approximate sparse nerves for classes of Dowker dissimilarities
  including all finite Dowker dissimilarities and Dowker
  dissimilarities whose homology is \v{C}ech persistent homology.

  All other sparsification methods and software packages that we are
  aware of calculate persistent homology with either an additive or a
  multiplicative interleaving. In \emph{\href{https://github.com/mbr085/Sparse-Dowker-Nerves}{dowker\textunderscore homology}}, we allow for any non-decreasing interleaving function
  \(\alpha\).

  We analyze the computational complexity of the algorithms and
  present some benchmarks. For Euclidean data in dimensions larger
  than three, the sizes of simplicial complexes we create are in
  general smaller than the ones created by SimBa. Especially when
  calculating persistent homology in higher homology dimensions, the
  differences can become substantial.

\end{abstract}

\section{Introduction} \label{sec:introduction}

Topological Data Analysis combines machine learning with topological
methods, most importantly persistent homology
\cite{Robins1999,Edelsbrunner2000}. The underlying idea is that data
has shape and this shape contains information about the
data-generating process \cite{Carlsson2009}. Persistent homology is a
method to characterize topological features that occur in data at
multiple scales. Its theoretical properties, in particular the
structure theorem and the stability theorem make persistent homology
an attractive machine learning method.

A major obstacle to the wide-spread use of persistent homology is its
computational complexity when analyzing large datasets. For example
the \v{C}ech complex grows exponentially with the number of points in
a point cloud. In order to be able to calculate persistent homology, a
number of approximations enable us to reduce the computational
complexity of persistent homology calculations
\cite{botnan15approximating,SRGeom,simba,DBLP:journals/corr/abs-1812-04966}.

Recently, Blaser and Brun have presented methods to sparsify nerves
that arise from general Dowker dissimilarities
\cite{SparseDowker,SFN}. In this article, we apply these techniques to
calculate the persistent homology of point clouds, weighted networks
and more general filtered covers. This paper is focused on the
algorithm implementation, computational complexity and benchmarking of
methods suggested in Blaser and Brun \cite{SFN}.

All algorithms presented in this manuscript are implemented in the
python package \emph{\href{https://github.com/mbr085/Sparse-Dowker-Nerves}{dowker\textunderscore homology}}, available on
\href{https://github.com/mbr085/Sparse-Dowker-Nerves}{github}.  With
\emph{\href{https://github.com/mbr085/Sparse-Dowker-Nerves}{dowker\textunderscore homology}} it is possible to calculate
persistent homology of ambient \v{C}ech filtrations, and intrinsic
\v{C}ech filtrations of point clouds, weighted networks and general
finite filtered covers. The \emph{\href{https://github.com/mbr085/Sparse-Dowker-Nerves}{dowker\textunderscore homology}}
package does all the preprocessing and sparsification, and relies on
\emph{\href{http://gudhi.gforge.inria.fr/}{GUDHI}} \cite{gudhi} for calculating persistent homology. Users
may specify additive interleaving, multiplicative interleaving or
arbitrary interleaving functions.

This paper is organized as follows. In Section~\ref{sec:theory}, we
give a short introduction on the underlying theory of the methods
presented here. Section~\ref{sec:algorithms} presents the implemented
algorithms in detail. In Section~\ref{sec:complexity} we quickly
discuss the size complexity of the sparse nerve and in
Section~\ref{sec:benchmarks} we provide detailed benchmarks comparing
the sparse Dowker nerve to other sparsification
strategies. Section~\ref{sec:conclusions} is a short summary of
results.

\section{Theory} \label{sec:theory}

The theory is described in detail in \cite{SFN}. In brief, the
algorithm consists of two steps, a truncation and a restriction. Given
a Dowker dissimilarity \(\Lambda\), the truncation gives a new Dowker
dissimilarity \(\Gamma\) that satisfies a desired interleaving
guarantee. The restriction constructs a filtered simplicial complex
that is homotopy equivalent to, but smaller than the filtered nerve of
\(\Gamma\).  The paper \cite{SFN} gives a detailed description of the
sufficient conditions for a truncation and restriction to satisfy a
given interleaving guarantee.  Here we give a new algorithm to choose
a truncation and restriction that together result in a small sparse
nerve. In Section~\ref{sec:benchmarks}, we compare sparse nerve sizes
from the algorithms presented here with the sparse nerve sizes of the
algorithms presented in \cite{SparseDowker} and \cite{SFN}.

\section{Algorithms} \label{sec:algorithms}

We present all algorithms given a finite Dowker
dissimilarity. Generating a finite Dowker dissimilarity from data is a
precomputing step that we do not cover in detail. For the intrinsic
\v{C}ech complex of \(n\) data points in Euclidean space
\(\mathbb{R}^d\), this consists of calculating the distance matrix,
with time complexity \(\mathcal{O}(n^2 \cdot d)\) operation.

\subsection{Cover matrix} \label{sec:covermatrix} 

The cover matrix is defined in \cite[Definition 5.4]{SFN}. Let
\(\Lambda \colon L \times W \to \rstar\) be a Dowker
dissimilarity. Given \(l,l'\in L\) let
\begin{displaymath}
  P(l,l') = \{ \Lambda(l',w) \, \mid \,
  w \in W \text{ with } \Lambda(l,w) < \Lambda(l',w) \}
\end{displaymath}
and define the cover matrix \(\rho\) as
\begin{displaymath}
  \rho(l,l') =
  \begin{cases}
    \sup P(l,l') & \text{if \(P(l,l')\) is non-empty} \\
    0 & \text{if \(P(l,l') = \emptyset\).} \\
  \end{cases}
\end{displaymath}

More generally, we can define a cover matrix of two Dowker
dissimilarities \(\Lambda_1 \colon L \times W \to \rstar\) and
\(\Lambda_2 \colon L \times W \to \rstar\) as follows.
\begin{displaymath}
  P(l,l') = \{ \Lambda_1(l',w) \, \mid \,
  w \in W \text{ with } \Lambda_2(l,w) < \Lambda_1(l',w) \}
\end{displaymath}
and define the cover matrix \(\rho\) as before. We define the cover
matrix algorithm in this generality, but sometimes we will use it with
just one Dowker dissimilarity \(\Lambda\), in which case we implicitly
use \(\Lambda_1 = \Lambda_2 = \Lambda\).

Our algorithms for calculating the truncated Dowker dissimilarity and
for calculating a parent function both rely on the cover matrix. The
cover matrix is the mechanism for the two algorithms to
interoperate. Algorithm~\ref{alg:covermatrix} explains how the cover
matrix can be calculated from two Dowker dissimilarities.

\begin{algorithm} \label{alg:covermatrix}
\SetKwInOut{Input}{Input}
\SetKwInOut{Output}{Output}
\Input{Dowker dissimilarities \(\Lambda_1(l, w)\) and
  \(\Lambda_2(l, w)\) for all \(l \in L\) and \(w \in W\).}
\Output{Cover matrix \(\rho(l_0, l_1)\) for all \(l_0, l_1 \in L\).}
\caption{Cover matrix}
Define \(\rho\) as an \(|L| \times |L|\) matrix of zeros indexed by
\(L \times L\). \\
\For{\((l_0, l_1)\) in \(L \times L\)}
{
  \For{\(w\) in \(W\)}
  {
  	\If{\(\Lambda_2(l_0, w) < \Lambda_1(l_1, w)\)}
  	{
  	  Update 
  	  \(\rho(l_0, l_1) = \max\{\rho(l_0, l_1), \Lambda_1(l_1, w)\}\). 
  	}
  }
}
Return \(\rho\).
\end{algorithm}

The cover matrix algorithm is the bottleneck for calculating the
truncated Dowker dissimilarity and the parent function. Its running
time \(\order(|L|^2 \cdot |W|)\) is quadratic in the size of \(L\) and
linear in the size of \(W\).

\subsection{Truncation} \label{sec:truncation}

Given a Dowker dissimilarity \(\Lambda \colon L \times W \to \rstar\),
and a translation function \(\alpha \colon \rstar \to \rstar\), every
Dowker dissimilarity \(\Gamma \colon L \times W \to \rstar\)
satisfying
\(\Lambda(l, w) \le \Gamma(l, w) \le \alpha(\Lambda(l, w))\), is
\(\alpha\)-interleaved with \(\Gamma\). In the case where \(\alpha\)
is multiplication by a constant, both extremes \(\Lambda(l, w)\) and
\(\alpha(\Lambda(l, w)\) will result in restrictions with sparse
nerves of the same size.  Our goal is to find a truncation that
interacts well with the restriction presented in Section
\ref{sec:restriction} in order to produce a small sparse nerve.

Algorithm~\ref{alg:truncation} explains in detail, how the truncated
Dowker dissimilarity is calculated. The high level view is that we
first calculate a farthest point sampling from the cover matrix and
the edge list \(E\) of the hierarchical tree of farthest points.
Finally, we iteratively reduce \(\Gamma(l, w)\) starting from
\(\alpha(\Lambda(l, w))\) by taking the minimum of \(\Gamma(l, w)\)
and \(\Gamma(l', w)\) for \((l', l)\) in \(E\).

\begin{algorithm} \label{alg:truncation}
\SetKwInOut{Input}{Input}
\SetKwInOut{Output}{Output}
\Input{Dowker dissimilarity \(\Lambda(l, w)\) for all \(l
  \in L\) and \(w \in W\), \\
\hspace{0pt}translation function \(\alpha \colon \rstar \to \rstar\). } 
\Output{Truncated dowker dissimilarity \(\Gamma(l, w)\) for all \(l
  \in L\) and \(w \in W\).} 
\caption{Truncated Dowker dissimilarity}
Calculate cover matrix \(\rho(l_0, l_1)\) of \(\Lambda\) and
\(\alpha \Lambda\) for all \(l_0, l_1 \in L\). \\
Choose initial point \(l_0 \in L\) and set \(L_0 = \{l_0\}\) and
\(T(l_0) = \infty\). \\
Initialize cover distance from \(L_0\) as \(d(L_0, l) = \rho(l,
l_0)\) for \(l \in L \setminus \{l_0\}\). \\
Set index \(i = 0\). \\
\While{\(|L_0| < |L|\)}
{
  Increment \(i\) by \(1\). \\
  Add the point \(l_{i} = \argmax_{l' \in L \setminus L_0} d(L_0, l')\) to
  \(L_0\). \\
  Set \(T(l_i) = d(L_0, l_i)\). \\
  Update the cover distance from \(L_0\) as \(d(L_0, l) =
  \min\{d(L_0\setminus \{l_{i}\}, l), \rho(l, l_{i})\}\) for \(l \in L
  \setminus L_0\). \\
}
Initialize the graph \(G = (L, E)\) with \(E = \emptyset\). \\
\For{\(l\) in \(L_0 \setminus \{l_0\}\) (sorted in order points were added to \(L_0\))} {
  \If{There exists a \(l' \in L\) with \(T(l) = \rho(l, l')\)}{
    Find the minimum \(\psi(l)\) such that \(T(l) = \rho(l, \psi(l))\). \\
  }
  \Else{
    Find the minimum of \(\rho(l, l')\) for \(l' < l\) in the order and the argument \(\psi(l)\) minimizing it. 
  }
  Add \((l, \psi(l))\) to the edge list \(E\). 
}
Topologically sort the nodes \(l \in L\) from highest to lowest \(T(l)\). \\
Initialize \(\Gamma(l, w) = \alpha(\Lambda(l, w))\) for \(l \in L\) and \(w \in W\). \\
\For{\(l\) in \(L \setminus \{l_0\}\) (topologically sorted)} {
  \For{\(l'\) such that \((l', l) \in E\)}{
  	Update \(\Gamma(l, -) = \min\{\Gamma(l, -), \Gamma(l', -)\}\). \\
  }
  Update \(\Gamma(l, -) = \max\{\Gamma(l, -), \Lambda(l, -)\}\). \\
}
Return \(\Gamma\).
\end{algorithm}

The truncation algorithm has a worst-case time-complexity
\(\order(|L|^2 \cdot |W|)\).  As mentioned earlier, calculating the
cover matrix is the bottleneck.  The time complexity of the
\emph{while} loop is \(\order(|L|^2)\), sorting is
\(\order(|L| \cdot \log |L|)\), the first \emph{for} loop is
\(\order(|L|^2)\), the topological sort of a tree is \(\order(|L|)\),
and the last \emph{for} loop is \(\order(|L| \cdot |W|)\).

\subsection{Parent function} \label{sec:parent}

The parent function \(\varphi \colon L \to L\) can in principle be any
function such that the graph \(G\) consisting of all edges
\((l, \varphi(l))\) with \(l \neq \varphi(l)\), is a tree.

Here we present the algorithm to create one particular parent function
that works well in practice and combined with the truncation presented
in Section \ref{sec:truncation} results in small sparse nerves.

Algorithm~\ref{alg:parent} is a greedy algorithm. Ideally, we would
like to set the parent point of any point \(l \in L\) as the point
\(l' \in L\) that minimizes \(\rho(l, l'')\) for \(l'' \in L\) with
\(\rho(l, l'') > 0\). However, this may not result in a proper parent
function. Therefore we start with this as a draft parent function and
then update it so that it becomes a proper parent function.

\begin{algorithm} \label{alg:parent}
\SetKwInOut{Input}{Input}
\SetKwInOut{Output}{Output}
\Input{Dowker dissimilarity \(\Lambda(l, w)\) for all \(l \in L\) and \(w \in W\).}
\Output{Parent points \(\varphi(l)\) for all \(l \in L\).}
\caption{Parent points}
Calculate cover matrix \(\rho(l_0, l_1)\) for all \(l_0, l_1 \in L\). \\
\For{\(l\) in \(L\)} {
  Find the minimum \(m(l)\) of \(\rho(l, l')\) for all \(l' \neq l\) 
  and the argument \(\varphi^{\ast}(l)\) which minimizes it.
}
Sort \(l \in L\) by non-increasing \(m(l)\). \\
Let \(l_0 \in L\) be the first point in \(L\). \\ 
Initialize \(\varphi(l) = l_0\) for all \(l \in L\). \\
\For{\(l\) in \(L  \setminus \{l_0\}\)} {
  \If{\(\varphi^{\ast}(l)\) comes before \(l\)}
  {
    Set \(\varphi(l) = \varphi^{\ast}(l)\).
  }
  \Else
  {
    Set \(\varphi(l) = \argmin \rho(l, l')\) for \(l'\) that come before
    \(l\) with \(\rho(l, l') > 0\).
  }
}
Return \(\varphi\).
\end{algorithm}

The time complexity of calculating the cover matrix is \(\order(|L|^2
\cdot |W|)\). Every subsequent step can be done in at most
\(\order(|L|^2)\) time. 

\subsection{Restriction} \label{sec:restriction}

Given a set of parent points \(\varphi(l)\) for \(l \in L\) and the
cover matrix \(\rho \colon L \times L \to \rstar\),
Algorithm~\ref{alg:restriction} calculates the minimal restriction
function \(R: L \to \rstar\) given in \cite[Definition 5.4,
Proposition 5.5]{SFN}.

\begin{algorithm} \label{alg:restriction}
\SetKwInOut{Input}{Input}
\SetKwInOut{Output}{Output}
\Input{Parent points \(\varphi(l)\) for all \(l \in L\), \\
\hspace{0pt}cover matrix \(\rho(l_0, l_1)\) for all \(l_0, l_1 \in L\). } 
\Output{Restriction times \(R(l)\) for all \(l \in L\).} 
\caption{Restriction times}
Initialize \(R'(l) = \infty\) for \(l \in L\). \\
\For{\(l\) in \(L\)}
{
  \If{\(\varphi(l)\) is not \(l\)}
  {
  Set \(R'(l) = \rho(l, \varphi(l))\).
  }
}

\For{\(l\) in \(L\)}
{
  Set \(R(l) = R'(l)\). \\
  Set \(l' = l\). \\ 
  \While{\(\varphi(l')\) is not \(l'\)}
  {
	Set \(l' = \varphi(l')\). \\ 
	Set \(R(l') = \max\{R(l'), R'(l')\}\).
  }
}
Return \(R\).
\end{algorithm}

The restriction algorithm has a worst-case quadratic time-complexity
\(\order(|L| ^ 2)\). The first loop is linear in the size of \(L\),
while the second loop depends on the depth \(td(G)\) of the parent
tree \(G\). For a given parent tree depth, the complexity is
\(\order(|L| \cdot td(G))\).

\subsection{Sparse Nerve} \label{sec:nerve}

In order to calculate persistent homology up to homological dimension
\(d\), we calculate the \((d+1)\)-skeleton \(N\) of the sparse
filtered nerve of   
\(\Gamma\). Given the truncated Dowker dissimilarity \(\Gamma\), the
parent tree \(\varphi\) and the restriction times \(R\),
Algorithm~\ref{alg:nerve} calculates the \((d+1)\)-skeleton
\(N\). Note that the filtration values
can be calculated either from \(\Gamma\) or directly from \(\Lambda\).

\begin{algorithm}[ht] \label{alg:nerve}
\SetKwInOut{Input}{Input}
\SetKwInOut{Output}{Output}
\Input{Dowker dissimilarities \(\Lambda(l, w)\) and \(\Gamma(l, w)\) for all \(l \in L\) and \(w \in W\), \\
\hspace{0pt}restriction times \(R(l)\) for all \(l \in L\), \\
\hspace{0pt}parent points \(\varphi(l)\) for all \(l \in L\), \\
\hspace{0pt}dimension \(d\)}
\Output{The \(d+1\)-skeleton \(N\) of the sparse nerve and filtration values \(v(\sigma)\) for \(\sigma \in N\).}
\caption{Sparse Nerve}
Initialize slope points \(S = L\). \\
\For{\(l\) in \(L\)}{
  Find the set \(L'\) of all points \(l' \in L\) with \(\varphi(l') = l\). \\
  Set \(r(L')\) to the maximum of \(R(l')\) for \(l' \in L'\). \\
  \If{\(R(l) < \infty\) and \(r(L') < R(l)\)}{
    Remove \(l\) from \(S\)
  }
}
Initialize maximal faces \(F\). \\
\For{\(l\) in \(L\)}{
  \For{\(w \in W\)}{
    \If{\(\Gamma(l, w) <= R(l)\)}{
      Find the face \(f\) consisting of all \(l' \in L\) with 
      \(R(l) \le R(l')\), \(\Gamma(l', w) \le R(l)\),  
      \(\Gamma(l', w) < \infty\), 
      and if \(l' \in S\), then 
      \(\Gamma(l', w) < R(l')\). \\
      Add \(f\) to \(F\). \\ 
    }
  }
}
Calculate the \(d+1\)-skeleton \(N\) of the sparse nerve consisting of all subsets \(\sigma\) of \(F\) of cardinality at most \(d+2\). \\
\For{\(\sigma\) in \(N\)} {
  Calculate the filtration value \(v(\sigma)\) of \(\sigma\) as
  \(v(\sigma) = \min_{w \in W}\max_{l \in \sigma} \Lambda(l, w)\).
}
Sort \(N\) by \(v(\sigma)\) for \(\sigma \in N\). \\
Return \(N\) and \(v\).
\end{algorithm}

The time complexity of the sparse nerve algorithm is
\(\order(|L|^2 \cdot |W| + |N|\log(|N|))\). The loop to find slope points had time complexity \(\order(|L|^2)\) The loop for finding maximal faces has a time
complexity of \(\order(|L|^2 \cdot |W|)\). The remaining operations
have time complexity \(\order(|N|\log(|N|)\).  Calculating persistent
homology using the standard algorithm is cubic in the number of
simplices.

So far we have considered the case of a Dowker dissimilarity
\(\Lambda \colon L \times W \to \rstar\) with finite \(L\) and
\(W\). This includes for example the intrinsic \v{C}ech complex of any
finite point cloud \(X\) in a metric space \((M, d)\), where
\(L = W = X\) and \(\Lambda = d\).

\subsection{Ambient \v{C}ech complex}

Let \(X\) be a finite subset of Euclidean space \(\mathbb{R}^n\) and
consider its ambient \v{C}ech complex.  For \(L = X\) and
\(W = \mathbb{R}^n\), the Dowker nerve of \(\Lambda = d|_{L\times W}\)
is the ambient \v{C}ech complex of \(X\).  Since \(W\) is not finite
we have to modify our approach slightly to in order to construct a sparse
approximation of the Dowker nerve of \(\Lambda\). 

We first calculate the restriction function \(R'(l)\) for \(l \in L\) of 
the intrinsic \v{C}ech complex \(\Lambda' = \Lambda|_{L\times L}\). Then 
we note that \(R(l) = 2R'(l)\) is a restriction function for
\(\Lambda\) \cite[Definition 5.3]{SFN}. We can use
Algorithm~\ref{alg:nerve} to calculate the   
simplicial complex \(N\) using the 
restriction times \(R\) and Dowker dissimilarity \(\Lambda'\). However, since \(W\) is infinite, we can not
directly compute the minimum used to calculate the filtration values
\(v(\sigma)\) for \(\sigma \in N\). 
We circumvent this problem by considering a filtered
simplicial complex \(K\) with the same underlying simplicial complex
as \(N\), but with filtration values inherited from the Dowker nerve
\(N\Lambda\).  This means that the filtration values are computed with
the miniball algorithm.  Thus, we construct a filtered simplicial
complex \(K\), such that, for all \(t \in [0, \infty]\) we have
\begin{displaymath}
  N_t \subseteq K_t \subseteq N\Lambda_t.
\end{displaymath}
Since \(N\) is \(\alpha\)-interleaved with \(N\Lambda\), it follows by
\cite[Lemma 2.14]{SFN} that
also \(K\) is \(\alpha\)-interleaved with \(N\Lambda\).

\subsection{Interleaving Lines}

\begin{figure}[ht]
\centering
\includegraphics[scale=0.8]{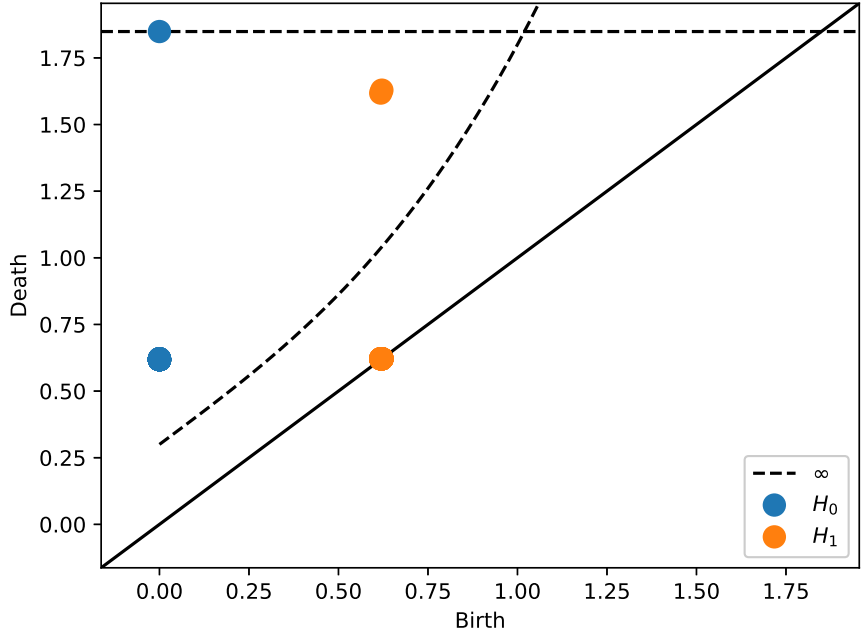}
\caption{Interleaving line. We generated \(100\) points on a Clifford
  torus that and calculated sparse persistent homology with an
  interleaving of \(\alpha(x) = \frac{x^3}{2} + x + 0.3\). This
  demonstrates the interleaving line for a general
  interleaving. Points above the line are guaranteed to have matching
  points in the persistence diagram with interleaving
  \(\alpha(x) = x\).}
\label{fig:interleaving}
\end{figure}

Our approximations to \v{C}ech- and Dowker nerves are interleaved with
the original \v{C}ech- and Dowker nerves. As a consequence their
persistence diagrams are interleaved with the persistence diagrams of
the original filtered complexes.  In order to visualize where the
points may lie in the original persistence diagrams, we can draw the
matching boxes from \cite[Theorem 3.9]{SFN}. However, this result in
messy graphics with lots of overlapping boxes. Instead of drawing
these matching boxes we draw a single interleaving line. Points
strictly above the line in the persistence diagram of the
approximation match points strictly above the diagonal in the
persistence diagram of the original filtered simplicial complex. More
precisely, the matching boxes of points above the interleaving line do
not cross the diagonal, while the matching boxes of all points below
the diagonal have a non-empty intersection with the
diagonal. Figure~\ref{fig:interleaving} illustrates such an
interleaving line for \(100\) data points on a Clifford torus with
interleaving \(\alpha(x) = \frac{x^3}{2} + x + 0.3\).

\section{Complexity Analysis} \label{sec:complexity}

We have shown time complexity analysis of each step. Combined, the
time it takes to calculate the sparse filtered nerve is
\(\order(|L|^2 \cdot |W| + |N|\log(|N|))\). Here we present some
results on the complexity of the nerve size depending on the maximal
homology dimension \(d\) and the sizes of the domain spaces \(L\) and
\(W\) of the Dowker dissimilarity \(\Lambda: L\times W \to \rstar\).
Although we can not show that the sparse filtered nerve is small in
the general case, we will show in the benchmarks below that this is
the case for many real-world datasets.

We now limit our analysis to Dowker dissimilarities that come from
doubling metrics and multiplicative interleavings with an interleaving
constant \(c>1\).  In that case, Blaser and Brun \cite{SFN} have showed 
that the size of the sparse nerve is bounded by the size of the
simplicial complex by Cavanna \textit{et al.} \cite{SRGeom}, whose
size is linear in the number \(|L|\) of points.

\section{Benchmarks} \label{sec:benchmarks}

We show benchmarks for two different types of datasets, namely data
from metric spaces and data from networks.  

\subsubsection{Metric data}
We have applied the presented algorithm to the datasets from Otter
\textit{et al.} \cite{Otter2017}.  First we split the data into two
groups, data in \(\mathbb{R}^d\) with dimension \(d\) at most \(10\)
and data of dimension \(d\) larger than \(10\).  The low-dimensional
datasets we studied consisted of six different Vicsek datasets
(Vic1-Vic6), dragon datasets with 1000 (drag1) and 2000 (drag2) points
and random normal data in 4 (rand4) and 8 (rand8) dimensions.  For all
low-dimensional datasets, we compared the sparsification method from
Cavanna et al. \cite{SRGeom} termed 'Sheehy', the method from
\cite{SparseDowker} termed 'Parent' and the algorithm presented in
this paper termed 'Dowker' for the intrinsic \v{C}ech complex. All methods were tested with a
multiplicative interleaving of \(3.0\).  In addition to the methods
described above, we have applied SimBa \cite{simba} with \(c = 1.1\)
to all datasets. Note that SimBa approximates the Rips complex with 
an interleaving guarantee 
larger than \(3.0\). For the \(3\)-dimensional data we additionally
compute the alpha-complex without any interleaving \cite{alpha}.  For
all algorithms we calculate the size of the simplicial complex used to
calculate persistent homology up to dimension \(1\)
(Table~\ref{tab:roadmap_benchmark_lowdim}).

\begin{table}
\centering
\begin{tabular}{lrrrrrrrr}
\hline
Name &  Points &  Dim &    Alpha &          Base &   Dowker &    Parent &     Sheehy &    SimBa \\
\hline
Vic1  &   300 &  3 &   5655 &  $4.5 \cdot 10^{6}$ &   1526 &   35371 &    29579 &    \textbf{830} \\
Vic2  &   300 &  3 &   5657 &  $4.5 \cdot 10^{6}$ &   1282 &   24977 &    25352 &    \textbf{812} \\
Vic3  &   300 &  3 &   5889 &  $4.5 \cdot 10^{6}$ &   1301 &   30894 &    27611 &    \textbf{822} \\
Vic4  &   300 &  3 &   5838 &  $4.5 \cdot 10^{6}$ &   1113 &   28722 &    24413 &    \textbf{804} \\
Vic5  &   300 &  3 &   5953 &  $4.5 \cdot 10^{6}$ &   1196 &   39098 &    68981 &    \textbf{973} \\
Vic6  &   300 &  3 &   6006 &  $4.5 \cdot 10^{6}$ &   1314 &   38860 &    67250 &    \textbf{971} \\
drag1 &  1000 &  3 &  21632 &  $1.7 \cdot 10^{8}$ &   6045 &  196660 &   201308 &   \textbf{3204} \\
drag2 &  2000 &  3 &  44446 &  $1.3 \cdot 10^{9}$ &  12230 &  534998 &   395740 &   \textbf{6368} \\
ran4  &   100 &  4 &        &  $1.7 \cdot 10^{5}$ &    \textbf{317} &    7356 &    36316 &    420 \\
ran8  &  1000 &  8 &        &  $1.7 \cdot 10^{8}$ &  \textbf{14126} &  598328 &  4366593 &  24980 \\
\hline
\end{tabular}
\label{tab:roadmap_benchmark_lowdim}
\caption{Comparison of sizes of simplicial complexes for homology
  dimension $1$ for low-dimensional datasets in Euclidean space. The
  smallest simplicial complexes in each dimension are displayed in
  bold. For all three-dimensional datasets, SimBa results in slightly
  smaller simplicial complexes. For the two datasets of dimensions 
  larger than three, the Dowker simplicial complex is smallest.}
\end{table}

The sparse Dowker nerve is always smaller than the sparse Parent and
sparse Sheehy nerves. In comparison to SimBa, it is noticeable that
the SimBa results in slightly smaller simplicial complexes if the data
dimension is three, but the sparse Dowker Nerve is smaller for most
datasets in dimensions larger than \(3\). For datasets of dimension
\(3\), the alpha complex without any interleaving is already smaller
than the Parent or Sheehy interleaving strategies, but Dowker
sparsification and SimBa can reduce sizes further.

The high-dimensional datasets we studied consisted of the H3N2 data
(H3N2), the HIV-1 data (HIV), the Celegans data (eleg), fractal
network data with distances between nodes given uniformly at random
(f-ran) or with a linear weight-degree correlations (f-lin), house
voting data (hou), human gene data (hum), collaboration network (net),
multivariate random normal data in 16 dimensions (ran16) and senate
voting data (sen).

For all high-dimensional datasets, we compared the intrinsic
\v{C}ech complex sparsified by the
algorithm presented in this paper ('Dowker') with a multiplicative
interleaving of \(3.0\) to the Rips complex sparsified by 
SimBa \cite{simba} with \(c = 1.1\). For
the high-dimensional datasets, we do not consider the 'Sheehy' and
'Parent' methods, because they take too long to compute and are
theoretically dominated by the 'Dowker' algorithm.  For all algorithms
we calculate the size of the simplicial complex used to calculate
persistent homology up to dimensions \(1\) and \(10\)
(Table~\ref{tab:roadmap_benchmark_highdim}).

\begin{table}
\centering
\begin{tabular}{lrr|rrr|rrr}
\hline
& & &  $1$-dimensional  & & & $10$-dimensional & & \\
Name &  Points &     Dim &    Base  &  Dowker  &     SimBa  &          Base  &   Dowker  &       SimBa  \\
\hline
H3N2   &  2722 &  1173 &  $3.4 \cdot 10^9$ &  \textbf{9478} &   11676 &  $3.4 \cdot 10^{32}$ &  \textbf{12503} &     25305 \\
HIV    &  1088 &   673 &  $2.1 \cdot 10^8$ &  \textbf{2972} &   14834 &  $5.5 \cdot 10^{27}$ &   \textbf{3273} &   1887483 \\
eleg   &   297 &   202 &  $4.4 \cdot 10^6$ &  \textbf{1747} &    2688 &  $8.2 \cdot 10^{20}$ &   \textbf{6229} &     14883 \\
f-lin  &   512 &   257 &  $2.2 \cdot 10^7$ &  \textbf{1651} &   10757 &  $6.1 \cdot 10^{23}$ &   \textbf{2927} &  13457079 \\
fr-ran &   512 &   259 &  $2.2 \cdot 10^7$ &  \textbf{1571} &   13419 &  $6.1 \cdot 10^{23}$ &   \textbf{2249} &  $\infty$ \\
hou    &   445 &   261 &  $1.5 \cdot 10^7$ &  \textbf{1168} &    2283 &  $1.1 \cdot 10^{23}$ &   \textbf{1233} &      3753 \\
hum    &  1397 &   688 &  $4.5 \cdot 10^8$ &  \textbf{4431} &  108118 &  $1.1 \cdot 10^{29}$ &   \textbf{5673} &  $\infty$ \\
net    &   379 &   300 &  $9.1 \cdot 10^6$ &  \textbf{1164} &    1207 &  $1.6 \cdot 10^{22}$ &   1617 &      \textbf{1425} \\
ran16  &    50 &    16 &  $2.1 \cdot 10^4$ &   \textbf{105} &     203 &  $1.7 \cdot 10^{11}$ &    \textbf{105} &       293 \\
sen    &   103 &    60 &  $1.8 \cdot 10^5$ &   \textbf{269} &     298 &  $1.8 \cdot 10^{15}$ &    \textbf{279} &       317 \\
\hline
\end{tabular}
\label{tab:roadmap_benchmark_highdim}
\caption{Comparison of sizes of simplicial complexes for homology
  dimensions $1$ and $10$ for high-dimensional datasets in Euclidean
  space. The smallest simplicial complexes in each dimension are
  displayed in bold. Except for one dataset, the Dowker
  sparsifications result in smaller simplicial complexes than
  SimBa. Note that we write $\infty$ when the computer ran out of
  memory.}
\end{table}

In comparison to SimBa, it is noticeable that the SimBa, the Dowker
Nerve is smaller for most datasets, with a more pronounced difference
for persistent homology in \(10\) dimensions.

\subsubsection{Graph data}

In order to treat data that does not come from a metric, we calculated
persistent homology from a Dowker filtration
\cite{2016arXiv160805432C}. Table~\ref{tab:graph_benchmark} shows the
sizes of simplicial complexes to calculate persistent homology in
dimensions \(1\) and \(10\) of several different graphs with \(100\)
nodes. In both cases we calculated persistent homology with a
multiplicative interleaving \(\alpha = 3\), and for the
\(1\)-dimensional case we also calculated exact persistent homology.
For the \(1\)-dimensional case, the base nerves are always of the same
size \(166750\), the restricted simplicial complexes for exact
persistent homology range from \(199\) to \(166750\), while the
simplicial complexes for interleaved persistent homology have sizes
between \(199\) and \(721\). The simplicial complexes to calculate
persistent homology in \(10\) dimensions do not grow much larger when
multiplicative interleaving is \(3\).

\begin{table}[ht]
\centering
\begin{tabular}{lrr|rrr|rr}
  \hline
  Data properties &&& \(1\)-d case &&& \(10\)-d case & \\
  \hline
  Name & Nodes & Edges & Base & Dowker & Dowker & Base & Dowker \\ 
  & & & & \(\alpha = 3.0\) & \(\alpha = 1.0\) & & \(\alpha = 3.0\) \\ 
  \hline
  Cycle graph & 100 & 100 & 166750 & 297 & 166750 & $1.2 \cdot 10^{15}$ & 305 \\ 
  Circular ladder graph &  & 150 &  & 324 & 166750 &  & 345 \\ 
  Ladder graph &  & 148 &  & 316 & 46894 &  & 333 \\ 
  Star graph &  & 99 &  & 199 & 199 &  & 199 \\ 
  Wheel graph &  & 198 &  & 199 & 199 &  & 199 \\ 
  Grid graph &  & 180 &  & 484 & 70286 &  & 721 \\ 
  Multipartite graph &  & 4000 &  & 199 & 166750 &  & 199 \\ 
  \((5\times20)\) &&&&& \\
   \hline
\end{tabular}
\caption{Comparison of sizes of simplicial complexes for homology
  dimensions $1$ and $10$ for graphs. For the \(1\)-dimensional case,
  we show that the Dowker restriction can in some cases reduce the
  simplicial complex significantly even without any truncation.}
\label{tab:graph_benchmark}
\end{table}

\section{Conclusions}
\label{sec:conclusions}

We have presented a new algorithm for constructing a sparse nerve and
have shown in benchmark examples that its size does not
grow substantially for increasing data or homology
dimension and that it in many cases outperforms SimBa.
In addition, the presented
algorithm is more flexible than previous sparsification strategies in
the sense that it works for arbitrary Dowker dissimilarities and
interleavings. We also provide a python package
\emph{\href{https://github.com/mbr085/Sparse-Dowker-Nerves}{dowker\textunderscore homology}} that implements the presented
sparsification strategy.

\paragraph{Acknowledgements}
This research was supported by the Research Council of Norway through
Grant 248840.

%
%

\begin{thebibliography}{8}
\providecommand{\url}[1]{\texttt{#1}}
\providecommand{\urlprefix}{URL }
\providecommand{\doi}[1]{https://doi.org/#1}

\bibitem{SparseDowker}
{Blaser}, N., {Brun}, M.: {Sparse Dowker Nerves}. ArXiv e-prints  (Feb 2018),
  \url{http://arxiv.org/abs/1802.03655}

\bibitem{SFN}
{Blaser}, N., {Brun}, M.: {Sparse Filtered Nerves}. ArXiv e-prints  (Oct 2018),
  \url{http://arxiv.org/abs/1810.02149}

\bibitem{botnan15approximating}
Botnan, M.B., Spreemann, G.: Approximating persistent homology in {E}uclidean
  space through collapses. Applicable Algebra in Engineering, Communication and
  Computing  \textbf{26}(1),  73--101 (2015). \doi{10.1007/s00200-014-0247-y}

\bibitem{Carlsson2009}
Carlsson, G.: Topology and data. Bull. Amer. Math. Soc. (N.S.)  \textbf{46}(2),
   255--308 (2009). \doi{10.1090/S0273-0979-09-01249-X}

\bibitem{SRGeom}
Cavanna, N.J., Jahanseir, M., Sheehy, D.R.: A geometric perspective on sparse
  filtrations. CoRR  \textbf{abs/1506.03797} (2015)

\bibitem{DBLP:journals/corr/abs-1812-04966}
Choudhary, A., Kerber, M., Raghvendra, S.: Improved topological approximations
  by digitization. CoRR  \textbf{abs/1812.04966} (2018).
  \doi{10.1137/1.9781611975482.166}

\bibitem{2016arXiv160805432C}
{Chowdhury}, S., {M{\'e}moli}, F.: {A functorial Dowker theorem and persistent
  homology of asymmetric networks}. arXiv e-prints arXiv:1608.05432 (Aug 2016)

\bibitem{simba}
Dey, T.K., Shi, D., Wang, Y.: Sim{B}a: an efficient tool for approximating
  {R}ips-filtration persistence via {\it sim}plicial {\it ba}tch-collapse. In:
  24th {A}nnual {E}uropean {S}ymposium on {A}lgorithms, LIPIcs. Leibniz Int.
  Proc. Inform., vol.~57, pp. Art. No. 35, 16 (2016).
  \doi{10.4230/LIPIcs.ESA.2016.35}

\bibitem{alpha}
Edelsbrunner, H., Kirkpatrick, D., Seidel, R.: On the shape of a set of points
  in the plane. IEEE Transactions on Information Theory  \textbf{29}(4),
  551--559 (July 1983). \doi{10.1109/TIT.1983.1056714}

\bibitem{Edelsbrunner2000}
Edelsbrunner, H., Letscher, D., Zomorodian, A.: Topological persistence and
  simplification. In: 41st {A}nnual {S}ymposium on {F}oundations of {C}omputer
  {S}cience ({R}edondo {B}each, {CA}, 2000), pp. 454--463. IEEE Comput. Soc.
  Press, Los Alamitos, CA (2000). \doi{10.1109/SFCS.2000.892133}

\bibitem{Otter2017}
Otter, N., Porter, M.A., Tillmann, U., Grindrod, P., Harrington, H.A.: A
  roadmap for the computation of persistent homology. EPJ Data Science
  \textbf{6}(1), ~17 (Aug 2017). \doi{10.1140/epjds/s13688-017-0109-5}

\bibitem{Robins1999}
Robins, V.: Towards computing homology from approximations. Topology
  Proceedings  \textbf{24} (01 1999)

\bibitem{gudhi}
{The GUDHI Project}: {GUDHI} User and Reference Manual. {GUDHI Editorial Board}
  (2015), \url{http://gudhi.gforge.inria.fr/doc/latest/}

\end{thebibliography}

%
\end{document}